\documentstyle{article}

\def\ni{\noindent }
\def\eq #1{(\ref{#1})}       
\def\l{\left}                   
\def\r{\right}                  
\def\fr{\frac}                  
\def\u{u}                       
\def\yx{$\{x\leftrightarrow y\}$}

\def\se #1{\S \ref{#1}}

\def\fz{f_0}
\def\fo{f_1}
\def\ft{f_2}
\def\f3{f_3}

\def\fl3{\tilde{\f3}}

\def\s3t{\tilde{{s_3}}}

\def\y1{\mbox{$y'$}}             

\begin{document}

\title{An Abel ordinary differential equation class generalizing known integrable classes}
\author{E.S. Cheb-Terrab$\,^{1,2}$,and A.D. Roche$^3$}

\date{}
\maketitle
\thispagestyle{empty}

\medskip
\centerline {$^1\,$ CECM, Department of Mathematics and Statistics,	Simon Fraser University} 
\centerline {$^2\,$ Department of Theoretical Physics, State University of Rio de Janeiro}
\centerline {$^3\,$Symbolic Computation Group, Faculty of Mathematics, University of Waterloo}

\bigskip
\centerline{\small (Received: 14 August 1999; revised 29 June 2000, 7 May 2002,  16 August 2002)}

\begin{abstract}

We present a multi-parameter non-constant-invariant 
class of Abel ordinary differential equations
with the following remarkable features. 
This one class is shown to unify, that is, contain as particular cases,
all the integrable classes presented by Abel, Liouville and Appell, as well
as all those shown in Kamke's book and various other references.
In addition, the class being presented includes other new and fully integrable subclasses,
as well as the most general parameterized
class of which we know whose members can systematically be mapped into
Riccati equations.
Finally, many integrable members of this class can be systematically mapped into an integrable member of
a different class. We thus find new integrable classes from
previously known ones.

\end{abstract}

\section{Introduction}

Abel-type differential equations of the first kind are of the form

\begin{equation}
\y1 = \f3\,  y^3 + \ft\,  y^2 + \fo\,  y + \fz,
\label{abel}
\end{equation}

\ni where $y\equiv y(x)$ and $f_i$ are
arbitrary functions of $x$. Abel equations appear in the reduction of order
of many second and higher order families \cite{kamke,green,sachdev}, and hence are frequently found in the
modelling of real problems in varied areas.
A general ``exact integration" strategy for these equations was first
formulated by Liouville \cite{liouville2} and is based on the
concepts of classes, invariants and the solution of the equivalence problem
\cite{appell,abel1,liouville3}. Generally speaking, two Abel equations of the first kind belong to
the same equivalence class if and only if one can be obtained from the other by
means of a transformation of the form

\begin{equation}
\{x = F(t),\ \ y= P(t)\,  \u + Q(t)\},
\label{tr}
\end{equation}

\ni where $t$ and $u\equiv u(t)$ are respectively the new independent and
dependent variables, and $F$, $P$ and $Q$ are arbitrary functions satisfying
$F'\neq0$ and $P\neq 0$. By changing variables $\{x=t,\ y=\left
(g_1\u+g_0\right )^{-1}\}$, where $\{g_1,g_0\}$ are arbitrary
functions of $t$, Abel equations of the first kind can always be written in
second-kind format

\begin{equation}
\y1 = \fr{\tilde{\f3}\, y^3 + \tilde{\ft}\, y^2 + \tilde{\fo}\, y +
\tilde{\fz}}{g_1\,y + g_0}.
\label{abel_2k}
\end{equation}

\ni Abel equations of the second kind belong to the same class as their first-kind
partners. However, due to the arbitrariness introduced when switching from
first to second kind, the transformation preserving the class for the
latter becomes

\begin{equation}
\{x = F(t),\ \ y= \fr{P_1(t)\,  \u + Q_1(t)}{P_2(t)\,  \u + Q_2(t)}\},
\label{tr_2k}
\end{equation}

\ni where $P_1Q_2-P_2Q_1 \neq 0$. There are infinitely many Abel equation
classes, and their classification is performed by means of algebraic
expressions invariant under \eq{tr}, the so called invariants, built with
the coefficients $f_i$ in \eq{abel} and their derivatives.



When the invariants of a given Abel equation are constant,
 its
integration is straightforward: the equation can be transformed into a 
separable equation as explained in textbooks \cite{murphy}. On the contrary,
when the invariants are non-constant, the integration strategy relies on
recognizing the equation as equivalent to one of a set of previously
known integrable equations, and then applying the equivalence transformation
to that equation. However, for non-constant-invariant Abel equations,
only a few integrable classes are known.
In a recent work \cite{abel1}, for instance, a classification of all
integrable cases presented by Abel, Liouville and others
\cite{abel,appell,halphen,liouville3}, including examples from Kamke's
book, showed, in all, only four classes depending on one parameter
and seven classes depending on no parameters.


In this work, we present a single multi-parameter Abel equation class
(AIA\footnote{The acronyms AIA, AIR and AIL are explained below.})
generalizing all the integrable classes collected in \cite{abel1}.
In addition, AIA contains as particular case a new subclass
(AIR), depending on 6 parameters, all of whose members can be
systematically transformed into Riccati-type equations. This AIR class, in
turn, includes a 4-parameter fully integrable subclass (AIL). Then we will
see in \se{AIA_theorem} that
many representatives of subclasses of this AIA class can be
mapped into non-trivial Abel equations belonging to different subclasses.
Hence, if the member being mapped is solvable, it can be used to generate
a different, maybe new, solvable and non-constant-invariant class as will be
shown in \se{fitting}.

Due to its simplicity and the potential preparation of computer algebra routines for
solving the related equivalence problem \cite{abel1}, the material
being presented appears to us as a convenient starting point for a more
systematic determination of exact solutions for Abel equations.

Another important differential equation problem, complementary to the one
discussed in this paper, is the one where a first-order equation $\y1 =
\Phi(x,y)$, where $\Phi$ is arbitrary, can be transformed into separable by
means of linear transformations \eq{tr}. This problem is a generalization of
the case where $\Phi$ is cubic in $y$ and the Abel equation has constant
invariant. An algorithmic approach for solving this problem is presented in
the precedent paper in this issue\cite{gts_ejam}.

\section{The AIL, AIR and AIA classes}
\label{GTI}

As mentioned in the introduction, the whole collection of integrable classes
presented in \cite{abel1}, consisting of four 1-parameter classes and 7
classes without parameters, can be obtained by assigning particular values
to the parameters of a single multi-parameter Abel class, AIA. In turn, AIA
contains a subclass, AIR, all of whose members can be transformed into
Riccati-type equations, and inside AIR there is a fully integrable subclass
(AIL) all of whose members can be mapped into first-order linear equations.

Since all these classes are obtained by means of the same procedure, to
better illustrate the ideas we discuss first this AIL (Abel, Inverse-Linear)
subclass. So, consider the general linear equation
$\y1+g(x)\,y+f(x)=0$, where $f$ and $g$ are arbitrary, after
changing variables by means of the {\it inverse} transformation 
\yx:\footnote{By \yx\ we mean changing variables $\{x = u,\,
y=t\}$ followed by renaming $\{u\rightarrow y,\, t\rightarrow x\}$.}

\begin{equation}
\y1=-\fr{1}{g(y)x+f(y)}.
\label{inverse_linear}
\end{equation}

\ni An implicit solution to this equation is easily expressed in terms of
quadratures as

\begin{equation}
C_1=x\,\exp\l({\int \!g(y)\ {dy}}\r)
+\int\!{\exp\l({\int \!g(y)\ {dy}}\r)}\,f(y)\ {dy}.
\label{inverse_linear_answer}
\end{equation}

\ni The key observation here is that \eq{inverse_linear} will be of
second-kind Abel type for many choices of $f$ and $g$. For instance, by taking

\begin{equation}
f(y)
 =
\fr {s_0 y+r_0} {a_3\, y^3+a_2\, y^2+a_1\, y+a_0},
\ \ \ \ \ \ \ \ 
g(y)
 = 
\fr {s_1 y+r_1} {a_3\, y^3+a_2\, y^2+a_1\, y+a_0},
\label{fg}
\end{equation}

\ni where $\{s_1,\,s_0,\,r_1,\,r_0,\,a_i\}$ are arbitrary constants, the resulting
Abel family is 

\begin{equation}
\y1=-\fr {a_3\, y^3+a_2\, y^2+a_1\, y+a_0}{\l (s_1\,x +s_0\r)\, y+r_1\,
x+r_0}.
\label{AIL_8}
\end{equation}

\ni This equation has non-constant invariant, and can be seen as a
representative of a non-trivial multi-parameter class, and from its
connection with a linear equation, its general solution is obtained directly
from Eqs.(\ref{inverse_linear_answer}) and (\ref{fg}). As shown in the
next subsection, four of the eight parameters in \eq{AIL_8} are superfluous,
and this number can be reduced further by splitting the class into
subclasses. Even so, the class is surprisingly large, including
multi-parameter subclasses we have not found elsewhere.\footnote{We noted
afterwards, however, that \eq{AIL_8} could also be obtained using a
different approach; for instance by following Olver \cite{olver} and
considering \eq{abel} as an ``inappropriate reduction" of a second-order
equation that has a solvable non-Abelian Lie algebra; the resulting
restrictions on the coefficients in \eq{abel} surprisingly lead to a
family of the same class as \eq{AIL_8}.} Among others, the class
represented by \eq{AIL_8} contains as particular cases the two integrable
1-parameter classes related to Abel's work \cite{abel,abel1}, and most of
the examples found in Kamke's as well as in other textbooks \cite{green}.

We note that if instead of starting with a linear equation we were to start
with a Bernoulli equation, instead of \eq{AIL_8} we would obtain

\begin{equation}
\y1=-{\frac {a_{{3}}{y}^{3}+a_{{2}}{y}^{2}+a_{{1}}y+a_{{0}}}{\left (s_1 x
+s_0 {x}^{\lambda}\right )y+ r_1x+ r_0 {x}^{\lambda}}},
\label{GTIB}
\end{equation}

\ni which is reducible to \eq{AIL_8} by changing
$\{x={t}^{1/(1-\lambda)},\,y=\u \}$ followed by a redefinition of the constants
$c_i \rightarrow (1-\lambda)\,a_i $, and so it belongs to the same class
as \eq{AIL_8}. 
However, if we start with a Riccati equation instead of a Bernoulli equation,
and hence instead of \eq{inverse_linear} we consider

\begin{equation}
\y1=-\fr{1}{h(y)\,{x}^{2}+g(y)\,x+f(y)},
\label{inverse_riccati}
\end{equation}

\ni and then choose $f(y)$, $g(y)$ and $h(y)$ as in \eq{fg},
we obtain a 10-parameter
$\{s_i,\,r_i,a_k\}$ Abel type family (AIR, meaning Abel, Inverse-Riccati)

\begin{equation}
\y1=-{\frac {a_{{3}}{y}^{3}+a_{{2}}{y}^{2}+a_{{1}}y+a_{{0}}}
{\left (s_2{x}^{2}+s_1 x+s_0\right )y+r_2{x}^{2}+r_1 x+r_0}},
\label{AIR_10}
\end{equation}

\ni which becomes a Riccati-type equation by changing \yx. Four of these ten
parameters are superfluous. The remaining 6-parameter class includes as
particular cases the parameterized families presented by Liouville and
Appell, as well as other classes depending on no parameters shown in Kamke
and \cite{green} and having solutions in terms of special functions. We note
that \eq{AIL_8} is also a particular case of \eq{AIR_10}.

Finally, the same ideas can be used to construct a more general Abel class
(AIA, meaning Abel, Inverse-Abel) embracing the previous classes
Eqs.(\ref{AIL_8}) and (\ref{AIR_10}) as particular cases. This AIA family
is obtained by taking as starting point an Abel equation of the second kind
\eq{abel_2k} and, taking the coefficients
$f_i$ of the form \eq{fg}, arriving at

\begin{equation}
\y1={\frac {\left (g_{{1}}a_{{3}}x+g_{{0}}a_{{3}}\right ){y}^{3}+
\left (g_{{1}}a_{{2}}x+g_{{0}}a_{{2}}\right ){y}^{2}+\left (g_{{1}}a_{
{1}}x+g_{{0}}a_{{1}}\right )y+g_{{1}}a_{{0}}x+g_{{0}}a_{{0}}}{\left (s
_{{3}}{x}^{3}+s_{{2}}{x}^{2}+s_{{1}}x+s_{{0}}\right )y+r_{{3}}{x}^{3}+
r_{{2}}{x}^{2}+xr_{{1}}+r_{{0}}}}.
\end{equation}

\ni Taking now $g_1$ and $g_0$ as constants and redefining $g_1 a_i
\rightarrow a_i$, $g_0 a_i \rightarrow b_i$, instead of \eq{AIR_10} we
obtain the 16-parameter Abel family

\begin{equation}
\y1=-{\frac {\left (a_{{3}}x+b_{{3}}\right ){y}^{3}+\left (a_{{2}}x+b_
{{2}}\right ){y}^{2}+\left (a_{{1}}x+b_{{1}}\right )y+a_{{0}}x+b_{{0}}
}{\left (s_{{3}}{x}^{3}+s_{{2}}{x}^{2}+s_{{1}}x+s_{{0}}\right
)y+r_{{3}}{x}^{3}+r_{{2}}x^2+r_{{1}}x+r_{{0}}}}.
\label{AIA_16}
\end{equation}

\ni Remarkably, by changing \yx, we
arrive at an equation of exactly the same type,

\begin{equation}
\y1=-{\frac {\left (s_{{3}}x+r_{{3}}\right ){y}^{3}+\left (s_{{2}}x+r_
{{2}}\right ){y}^{2}+\left (s_{{1}}x+r_{{1}}\right )y+s_{{0}}x+r_{{0}}
}{\left (a_{{3}}{x}^{3}+a_{{2}}{x}^{2}+a_{{1}}x+a_{{0}}\right )y
+b_{{3}}{x}^{3}+b_{{2}}{x}^{2}+b_{{1}}x+b_{{0}}}}.
\label{inv_AIA}
\end{equation}

\ni Although we are not aware of a method for solving this AIA class for
arbitrary values of its 16 parameters, by applying the change of variables
\yx\ to representatives of solvable subclasses of AIA, one may obtain
representatives of non-trivial new solvable Abel classes (see
\se{AIA_theorem}). Due to this feature and since AIA already contains as a
particular case the AIR class and hence AIL (Eqs.(\ref{AIR_10}) and
(\ref{AIL_8})), this AIA class generalizes in one all the solvable classes
presented in \cite{abel1}, collected from the literature (see
\se{fitting}).

\subsection{The intrinsic parameter-dependence of the classes AIL and AIR}


For the purposes of finding exact solutions to Abel equations, and in so doing,
solving the equivalence problem with respect to the classes
represented by (\ref{AIL_8}) and to some extent (\ref{AIR_10}), it is
relevant to determine on how many parameters these classes
intrinsically depend. In fact, the technique we have been using for
tackling the equivalence problem relies heavily on the calculation of
multivariate resultants, and even for classes depending on just one
parameter these calculations require the use of special techniques 
\cite{abel1}.


Although the procedure presented below does not exhaust the possible reduction of
the number of parameters found in (\ref{AIL_8}) and (\ref{AIR_10}), it
suffices to reduce this number by four.
The idea is to explore fractional linear transformations of the form

\begin{equation}
\{x=t,\  y(x)=\frac{a\,u + b}{c\,u+d} \}.
\end{equation}

\ni We note that this transformation does not change the degrees in $x$ of
the denominators of the right-hand sides of (\ref{AIL_8}) and
(\ref{AIR_10}). Keeping that fact in mind, we rewrite the representatives
of the AIL and AIR classes as

\begin{equation}
\y1={\frac {\left (y-\alpha_{{0}}\right )\left (y-\alpha_{{1}}\right )
\left (y-\alpha_{{2}}\right )}{G_{{1}}(x)\,y+G_{{0}}(x)}},
\label{moebius}
\end{equation}

\ni where $\{G_{{1}}(x),\, G_{{0}}(x)\}$ are polynomials in
$x$ of the same degree (the ones implied by either \eq{AIL_8} or \eq{AIR_10}). By changing variables $\left
\{x=t,\, y=u(t)+\alpha_{{0}}\right \}$, \eq{moebius} becomes

\begin{equation}
\y1={\frac {y\left (y-\Delta_{{1,0}}\right )\left (y-\Delta_{{2,0}}
\right )}{G_{{1}}(x)\,y+\tilde{G_0}(x)}},
\label{seed_delta}
\end{equation}

\ni where $\Delta_{{i,0}} \equiv \alpha_i - \alpha_0$ and $\tilde{G_0}(x)$
is a polynomial in $x$ of the same degree as $G_0(x)$. Three cases now
arise, depending whether three, two or none of the roots $\alpha_i$
are equal.

\smallskip
\ni {\underline{\it Case 1: $\alpha_0=\alpha_1=\alpha_2$}}
\smallskip

In such a case, \eq{moebius} already depends on four fewer parameters than
\eq{AIL_8} or \eq{AIR_10}; a further change of variables
$\{x=t,y=1/u\}$ transforms \eq{moebius} to the form

\begin{equation}
\y1={\frac {1}{M_{{1}}(x)\,y+{M_0}(x)}},
\label{seed_delta_1}
\end{equation}

\ni where $\{M_{{1}}(x),\,M_{{0}}(x)\}$ are polynomials in $x$ of the same
degree as $\{G_{{1}}(x),\,G_{{0}}(x)\}$.

\smallskip
\ni {\underline{\it Case 2: $\alpha_0 = \alpha_1 \neq \alpha_2$}}
\smallskip

In this case, $\Delta_{{1,0}}=0$, so that changing variables
$$
\{x=t,y={\frac {1}{u(t)+{\Delta_{{2,0}}}^{-1}}} \}
$$
in \eq{seed_delta} leads to an equation of the form

\begin{equation}
\y1={\frac {y}{\tilde{M_{{1}}}(x)\,y+{\tilde{M_0}}(x)}},
\label{seed_delta_2}
\end{equation}

\ni where $\{\tilde{M_{{1}}}(x),\,\tilde{M_{{0}}}(x)\}$ are polynomials in
$x$ of the same degree as $\{G_{{1}}(x),\,G_{{0}}(x)\}$. Note that \eq{seed_delta_2}
also depends on four fewer parameters than \eq{AIL_8} or
\eq{AIR_10}.

\smallskip
\ni {\underline{\it Case 3: $\alpha_0 \neq \alpha_1 \neq \alpha_2$}}
\smallskip

In this case, $\Delta_{{1,0}}\neq \Delta_{{2,0}} \neq 0$, and by changing variables
$$
\{x=t,\ y=
\frac{1}{\left ({\Delta_{{2,0}}}^{-1}-{\Delta_{{1,0}}}^{-1}\right )u(t)+{\Delta_{{1,0}}}^{-1}}
\}
$$
in \eq{seed_delta} we obtain

\begin{equation}
\y1={\frac {y\,(y-1)}{\tilde{\tilde{M_{{1}}}}(x)+{\tilde{\tilde{M_0}}(x)}\,y}},
\label{seed_delta_3}
\end{equation}

\ni where $\{\tilde{\tilde{M_{{1}}}}(x),\,\tilde{\tilde{M_{{0}}}}(x)\}$
are polynomials in $x$ of the same degree as $\{G_{{1}}(x),\,G_{{0}}(x)\}$
and so \eq{seed_delta_3} too depends on four fewer parameters than either
\eq{AIL_8} or \eq{AIR_10}.


\subsection{Splitting of AIL into cases}
\label{splitting}

Since the AIL class is fully solvable, it makes sense, for the purpose of
tackling its related equivalence problem, to consider the further maximal
reduction in the number of parameters, and hence completely split the
class into a set of non-intersecting subclasses, all of which depend
intrinsically on a minimal number of parameters. With this motivation, we
performed some algebraic manipulations, finally determining that
the AIL class actually consists of two different subclasses,
respectively depending on two and one parameters.


To arrive at convenient representatives for these two subclasses of AIL, we
start by changing variables\footnote{We are assuming $\omega \equiv
r_1s_0-r_0s_1 \neq 0$ and $s_1 \neq 0$, which is justified, since when
$\omega = 0$ \eq{AIL_8} has constant invariant, therefore presenting
no interest, and when $s_1 = 0$ it can be transformed into an equation of the
form \eq{AIL_4} anyway by means of $\{x=t/r_1, y = (u-r_0)/s_0\}$.} in
\eq{AIL_8} according to

\begin{equation}
\{
x=-{\frac {t}{{s_1}^{2}}},\ 
y=-{\frac {r_0{s_1}^{2}+r_1\u }{\left (s_1s_0+\u \right )s_1}}\},
\label{tr_AIL_8_4}
\end{equation}

\ni and introducing new parameters $\{k_0,\,k_1,\,k_2,\,k_3\}$ according to

\begin{eqnarray}
{k_0} & = &{\frac
{\left (
    \left ({a_2}\,{r_0}^{2} + \left (s_0{a_0}-{a_1}\,r_0\right )s_0\right )s_0
    - {a_3}\,{r_0}^{3}
\right ){s_1}^{2}
}
{\left (r_1s_0-r_0s_1\right )^{2}}},
\\*[.13 in]
{k_1} & = & {\frac {
\left (
    {a_2}\,s_1 - 3\,{a_3}\,r_1\right ){r_0}^{2}
    + \left (
        \left (2\,r_1{a_2} - 2\,{a_1}\,s_1\right )r_0
        + \left (3\,{a_0}\,s_1- {a_1}\,r_1\right )s_0
      \right )s_0}
{\left (r_1s_0-r_0s_1
\right )^{2}}},
\\*[.13 in]
{k_2} & = &
{\frac
    {
    \left (
        \left (2\,{a_2}\,s_1 - 3\,{a_3}\,r_1 \right )r_1
        - {a_1}\,{s_1}^{2}
    \right)r_0
    + \left (3\,{s_1}^{2}{a_0} + \left (r_1{a_2}- 2\,{a_1}\,s_1\right )r_1\right )s_0
    }
{{s_1}^{2
}\left (r_1s_0-r_0s_1\right )^{2}}},
\\*[.13 in]
{k_3} & = & {\frac {
    {a_0}\,{s_1}^{3}
    + \left ( 
      \left ({a_2}\,s_1 - {a_3}\,r_1\right )r_1
      - {a_1}\,{s_1}^{2}
      \right )r_1}
{{s_{{1}}}^{4}\left (r_{{1}}s_{{0}}-r_{{0}}s_{{1}}\right )^{2}}.       
}
\end{eqnarray}

\ni Thus \eq{AIL_8} is transformed into

\begin{equation}
\y1={\frac {k_{{3}}{y}^{3}+k_{{2}}{y}^{2}+k_{{1}}y+k_{{0}}}{y+x}},
\label{AIL_4}
\end{equation}

\ni whose solution can be obtained directly from
Eqs.(\ref{inverse_linear_answer}) and (\ref{fg}) by taking $s_1=0,\,
r_1=1,\, s_0=1,\, r_0=0$ and $a_i \rightarrow -k_i$. For classification
purposes (see \se{fitting}), it is convenient to write \eq{AIL_4} in first-kind
form by changing variables $\{x=t,y=\frac{1}{u}-t\}$, leading to

\begin{equation}
\y1=
(k_{{3}}{x}^{3} - k_{{2}}{x}^{2} + k_{{1}}x -k_{{0}})\, {y}^{3}
- (3k_{{3}}{x}^{2} - 2k_{{2}}x + k_{{1}} + 1)\, {y}^{2}
+(3k_{{3}}x-k_{{2}})\,y-k_{{3}}.
\label{AIL_FirstKind}
\end{equation}

\ni Two different cases now arise, leading to two non-intersecting subclasses.

\smallskip
\ni {\underline{\it Case $k_3 \neq 0$}}
\smallskip

By redefining $k_3 \equiv -{k_4}^2$, then changing
variables in \eq{AIL_4} according to

\begin{equation}
\{x = -\fr{k_2 + 3\, t\, k_4}{3\, {k_4}^2},\  y=-\fr{1}{k_4\, \u}+\fr{k_2+3\, t\, k_4}{3\, {k_4}^2}\},
\label{composed_tr}
\end{equation}

\ni and next redefining $\{k_4,k_0,k_1\}$ in terms of new parameters
$\{\alpha,\,\beta,\,\gamma\}$ according to

\begin{equation}
k_4 = -\fr{\beta}{\gamma},\ \ \ \ \ \ \ \ \ 
k_0 = \fr{{k_2}^3 \gamma^4}{27 \beta^4}-\fr{k_2 \alpha \gamma^2}{3 \beta^2}+\gamma,
\ \ \ \ \ \ \ \ \ 
\label{beta_eq}
k_1 = \alpha - \fr{{k_2}^2 \gamma^2}{3 \beta^2},
\nonumber
\end{equation}


\ni we obtain a 2-parameter representative of the class, already in first kind format:

\begin{equation}
\y1=\left (x\alpha-\beta-{x}^{3}\right ){y}^{3}+\left
(3\,{x}^{2}-1-\alpha\right ){y}^{2}-3\,xy+1.
\label{AIL_2}
\end{equation}

\smallskip
\ni {\underline{\it Case $k_3 = 0$}}
\smallskip

This other branch of \eq{AIL_4} splits into two subcases: $k_2=0$ and $k_2
\neq 0$. In the former case \eq{AIL_4} becomes constant-invariant,
and is thus of no interest. When $k_2 \neq 0$, by introducing a new
parameter $\alpha$ by means of

\begin{equation}
\alpha=k_{{2}}k_{{0}}-\fr {{k_{{1}}}^{2}}{4},
\end{equation}

\ni and changing variables in \eq{AIL_4} according to

\begin{equation}
\{
x={\frac {k_{{1}}-2\,t}{2\,k_{{2}}}},\ 
y={\frac {2\,t-k_{{1}}}{2\,k_{{2}}}}-{\frac {1}{k_{{2}}\u }}
\},
\end{equation}

\ni one obtains a simpler representative for this class, depending on just
one parameter $\alpha$:

\begin{equation}
\y1=\left (\alpha+{x}^{2}\right ){y}^{3}
-\left (2\,x+1\right ){y}^{2}+y.
\label{AIL_1}
\end{equation}

\subsection{Generating new integrable classes from solvable members of AIA}

\label{AIA_theorem}

The motivation for this work was to try to 
unify the integrable Abel classes we have seen in the
literature. This goal has been partially accomplished with the formulation of the
AIL and AIR Abel classes, but there are still other integrable
classes, not included in AIR, which however all have the following property:
these Abel classes have representatives which could be obtained by changing variables \yx\
in the representative of other Abel classes. In this sense, these representatives are both of
Abel and inverse-Abel types, which led us to ask the following
question:

\begin{quote}
{\it Which Abel classes lead to other Abel classes by applying the
inverse transformation \yx\ to one of its representatives?}
\end{quote}

\ni We have formulated the answer to this question in terms of the following proposition
and its corollary.

\medskip \ni {\bf Proposition}: {\em If one Abel equation, $\alpha$, maps
into another Abel equation, $\beta$, by means of the inverse transformation
\yx, then both are of the form \eq{AIA_16}.} \medskip

\medskip \ni {\bf Proof}: By hypothesis $\beta$ is both of Abel and
inverse-Abel type, that is, it is of the form $\y1 = G(x,y)$, where $G$ is
both cubic over linear in $y$ (the Abel condition) and linear over cubic in
$x$ (the inverse-Abel condition). Hence, $G$ is a rational function of $x$
and $y$, with numerator cubic in $y$ and linear in $x$, and denominator
cubic in $x$ and linear in $y$, and so $\beta$ is in fact of the form
\eq{AIA_16}.

\medskip \ni {\bf Corollary}: {\em If a given Abel class is related
through the inverse transformation to another Abel class then both classes
have representatives of the form \eq{AIA_16}.} \medskip 

One consequence of this proposition is that AIA, \eq{AIA_16}, is in fact the
most general equation which is both Abel and inverse-Abel; similarly
AIR, \eq{AIR_10}, is the most general which is both Abel and
inverse-Riccati; and AIL, \eq{AIL_8}, is the most general which is both Abel
and inverse-linear.

It is worth mentioning here that an Abel class may have many different
representatives of the form \eq{AIA_16}. Consequently, for instance, the
AIL integrable class represented by \eq{AIL_4}, which naturally maps into a
linear equation by means of \yx, also contains members which map into
non-trivial Abel equations by means of the same transformation. As an example
of this, consider 151 of Kamke's book:

\begin{equation}
\y1={\frac {1-2\,xy+{y}^{2}-2\,{y}^{3}x}{{x}^{2}+1}}.
\label{this_151}
\end{equation}

\ni This equation has the form \eq{AIA_16} and, by changing variables \yx, it
leads to a non-trivial new Abel class and nevertheless it belongs to the AIL
class (see \eq{class_7}). In other words, the class represented by this equation has 
a representative of the form \eq{AIA_16} different from \eq{this_151},
which, by means of \yx, maps into a linear equation.
An example where the same thing happens with a member of the AIR class is
given by

\begin{equation}
\y1=-2\,x{y}^{2}+{y}^{3}.
\label{this_one}
\end{equation}

\ni This equation was presented by Liouville \cite{liouville2}; it is of the form
\eq{AIA_16} and by means of \yx\ leads to another non-trivial Abel class.
In addition \eq{this_one} belongs to the AIR class, so that it has another
representative of the form \eq{AIA_16}, which by means of \yx\ maps into a
Riccati equation. Some other examples illustrating how new solvable classes can
be obtained from solvable members of AIA by changing variables \yx\ are
shown in the next section.

\section{AIA: a generalization of known integrable classes}
\label{fitting}

As mentioned in the Introduction, all the solvable classes collected in
\cite{abel1}, four depending on one parameter, labelled A, B, C and D, and
another seven not depending on parameters, labelled 1 to 7, are particular
members of the class represented by \eq{AIA_16}. The ``classification" of
these solvable classes as particular members of AIL, AIR or AIA is as
follows:

\begin{table}
\caption{Classification of the integrable classes collected in \cite{abel1}.}
\begin{tabular}{                     
cccc
}
\hline
\hline
 & Subclass AIL & Subclass AIR & Class AIA  \\
\hline
Classes & A, C, 4, 5 & B, D, 2 &  1, 3, 6, 7 \\
\hline
\hline
\end{tabular}
\end{table}

Although this classification shown in Table 1 is not difficult to verify,
it is worthwhile to show it explicitly. Starting with the parameterized
class by Abel \cite{abel}, shown in \cite{abel1} as ``Class A", 

\begin{equation}
\y1=\left (\alpha\,x+\fr{1}{x}+\fr{1}{x^{3}}\right ){y}^{3}+{y}^{2},
\label{class_A}
\end{equation}

\ni where $\alpha$ is the parameter, this equation can be obtained from
\eq{AIL_FirstKind} by taking
$
\{k_{{3}}=2\,\alpha,k_{{2}}=-1,k_{{1}}=1/2,k_{{0}}=0 \}
$
and changing
variables
$$
\{x=\fr{{t}^{2}}{2},\ y=2\,{\frac {\u +t}{{t}^{3}}}\}.
$$

\ni So, \eq{class_A} is a member of AIL.

Concerning the 1-parameter class by Liouville \cite{liouville3}, labelled
in \cite{abel1} as class ``B",

\begin{equation}
\y1=2\,\left ({x}^{2}-\alpha\right ){y}^{3}+2\,\left (x+1\right ){y}^{2},
\label{class_B}
\end{equation}

\ni this equation is obtained from \eq{AIR_10} by taking
$\{s_{{2}}=0,s_{{1}}=0,s_{{0}}=1,r_{{2}}=1,r_{{1}}=0,r_{{0}}=0,a_{{3}}=0,
a_{{2}}=0,a_{{0}}=-2\,\alpha,a_{{1}}=-2\}$ and changing variables

$$
\{x=t,\ y=-\fr{1}{\u }-{t}^{2} \}.
$$

\ni \eq{class_B} is then a member of AIR.

The next solvable class, related to Abel's work, presented in \cite{abel1}
as class ``C", is represented by

\begin{equation}
\y1={\frac {\alpha\left (1-{x}^{2}\right ){y}^{3}}{2\,x}}+\left (\alpha-1
\right ){y}^{2}-{\frac {\alpha y}{2\,x}},
\label{class_C}
\end{equation}

\ni and this equation can be obtained from \eq{AIL_FirstKind} by taking
$\{k_{{2}}=0,\,k_{{1}}=\alpha/2,\,k_{{0}}=0\,\,k_{{3}}=-1/2\}$ and changing variables
$$
\{x={\frac {\sqrt {\alpha}}{t}},\ y={ \frac{t\left (1-t\,\u \right)}{\sqrt{\alpha}}} \}.
$$
\eq{class_C} is then a member of AIL.

The last parameterized solvable class shown in \cite{abel1}, labelled
there as class ``D", is related to Appell's work \cite{appell}, and is given
by

\begin{equation}
\y1=-{\frac {{y}^{3}}{x}}-{\frac {\left (\alpha+{x}^{2}\right ){y}^{2}}{{x}
^{2}}}.
\label{class_D}
\end{equation}

\ni This equation is obtained from \eq{AIR_10} by taking $\{
s_{{2}}=0,s_{{1}}=1,s_{{0}}=0,r_{{2}}=1,r_{{1}}=0,r_{{0}}=-\alpha,a_{{3}}=0,
a_{{2}}=0,a_{{1}}=0,a_{{0}}=-1\}$ and changing variables
$$
\{x=t,\ y=\fr{1}{\u }+{\frac {\alpha-{t}^{2}}{t}} \}.
$$
So, \eq{class_D} is a member of AIR.

Concerning the classes collected in \cite{abel1} not depending on
parameters, the first one, there labelled as ``Class 1", was presented by
Halphen \cite{halphen} in connection with doubly periodic elliptic
functions:

\begin{equation}
\y1={\frac {3\,y\left (1+y\right )-4\,x}{x\left (8\,y-1\right )}}.
\label{class_1}
\end{equation}

\ni This equation, clearly a member of \eq{AIA_16}, can be obtained by changing
\yx\ in
$$
\y1={\frac {y\left (8\,x-1\right )}{3\,x\left(x+1\right )-4\,y}}
$$
\ni which in turn can be obtained from \eq{class_C} (the
solvable AIL class) by taking $\alpha=- 2/3$ and changing variables
$$
\{x =\fr{(1-2\,t)\sqrt{3-6\,t}}{9\,t},\
y = {\frac {12\,ut\sqrt {3-6\,t}}{\left (t+1\right )^{2}\left(4\,u - 3\,{t}^{2} - 3\,t\right )}}
 \}.
$$
\ni This derivation also illustrates how
new solvable classes can be obtained by interchanging the roles between
dependent and independent variables in solvable members of AIA.

As for the representative of Class 2, by Liouville \cite{liouville2},
shown in \cite{abel1} as

\begin{equation}
y'=y^3-2\,x\,y^2,
\label{class_2}
\end{equation}
this equation is obtained from \eq{AIR_10} by taking
$\{s_{{2}}=0,s_{{1}}=0,s_{{0}}=1,r_{{2}}=1,r_{{1}}=0,r_{{0}}=0,a_{{3}}=0,
a_{{2}}=0,a_{{1}}=0,a_{{0}}=1\}$ and converting the resulting equation into
first-kind format by changing variables
$$
\{x=t,\ y=\fr{1}{\u }-{t}^{2}\}.
$$
\ni \eq{class_2} is then a member of AIR.

Also of note here, \eq{class_2} is a special case of the equation presented by
Appell in \cite{appell}:

\begin{equation}
\y1= -\fr{y^3}{\alpha\,x^2+\beta\,x+\gamma} 
-\fr{d}{dx}\l(\fr{ax^2+bx+c}{\alpha\,x^2+\beta\,x+\gamma}\r) y^2
\label{general_appell}
\end{equation}

\ni with $\alpha=0,\,\beta=0,\,\gamma=-1,\,a=-1,\,b=0,\,c=0$. \eq{general_appell} is also
seen to be a member of AIR since it can be obtained by changing
$$
\{x=t,\ y=\fr{1}{u}-\fr{a\,t^2+b\,t+c}{\alpha\,t^2+\beta\,t+\gamma}\}
$$
\ni in \eq{AIR_10} with
$s_2=\alpha,\,s_1=\beta,\,s_0=\gamma,\,r_2=a,\,r_1=b,\,r_0=c,\,a_3=0,\,a_2=0,\,a_1=0,\,a_0=-1$.

Class 3, also by Liouville \cite{liouville3}, presented in \cite{abel1} as

\begin{equation}
\y1={\frac {{y}^{3}}{4\,{x}^{2}}}-{y}^{2}
\label{class_3}
\end{equation}

\ni can be obtained from \eq{class_2} by changing \yx\ and then converting
it to first-kind format by means of
$$
\{x=2\,t,\ y=-\fr{1}{\u }+t \}.
$$
\ni \eq{class_3} is then a member of AIA. This also illustrates the derivation of
a solvable class by changing \yx\ in solvable members of the AIR
subclass.

The next class, presented in
\cite{abel1} as Class 4, collected among the Abel equation examples of Kamke's book,

\begin{equation}
\y1={y}^{3}-{\frac {\left (x+1\right ){y}^{2}}{x}},
\label{class_4}
\end{equation}

\ni can be obtained from \eq{AIL_1} by taking $\alpha=0$ and changing
variables
$$
\{x=-\fr{1}{t},\ y={t}^{2}\u -t\},
$$
\ni so it belongs to the AIL subclass.

In \cite{abel1}, three new integrable classes not depending on
parameters were presented too - these are classes ``5", ``6" and ``7".
Starting with class 5, given by

\begin{equation}
\y1=-{\frac {\left (2\,x+3\right )\left (x+1\right ){y}^{3}}
{2\,{x}^{5}}}+{\frac {\left (5\,x+8\right ){y}^{2}}{2\,{x}^{3}}},
\end{equation}

\ni this equation can be obtained from \eq{AIL_8} by taking
$\{s_{{1}}=0,\,s_{{0}}=1,\,r_{{1}}=1,\,r_{{0}}=0,\,a_{{3}}=-6,\,a_{{2}}=10,\,a_{{1}}=-4,\,
a_{{0}}=0\}$ and changing variables
$$
\{x=\fr{1}{t},\ 
y={\frac {2\,{t}^{2}-\left(t+1\right )\u}{t \left (\left (t+1\right )\u+2\,t\right
)}} \}.
$$

Regarding Class 6, given by

\begin{equation}
\y1=-{\frac {{y}^{3}}{{x}^{2}\left (x-1\right )^{2}}}+{\frac {\left (
1-x-{x}^{2}
\right ){y}^{2}}{{x}^{2}\left (x-1\right )^{2}}},
\label{class_6}
\end{equation}
this equation is obtained from \eq{class_4} by changing \yx\ and then converting
it to first-kind form by means of
$$
\{x=t,\,y={\frac {t-1+\u }{\left(t-1\right )\u }}\}.
$$

Concluding, class 7 is given by

\begin{equation}
%
\y1={\frac {\left (4\,{x}^{4}+5\,{x}^{2}+1\right ){y}^{3}}{2\,{x}^{3
}}}+{y}^{2}
+{\frac {\left (1- 4\,{x}^{2}\right )y}{2\,x\left ({x}^{2}
+1\right )}}.
\label{class_7}
\end{equation}

\ni This equation can be obtained from Kamke's first order example 151,
\begin{equation}
\y1={\frac {\left ({y}^{2}+1\right )\left (1-2\,yx\right )}{{x}^{2}+1
}},
\label{k151}
\end{equation}

\ni by first changing \yx\ and then converting
the resulting equation to first kind format by means of
$$
\{x=t,\ y=\fr{1}{2\,t} - \fr{1}{2\,\l(t^2+1\r)\, u}\}.
$$

\ni In turn, \eq{k151} is a member of the
AIL subclass and can be obtained from \eq{class_C} by taking $\alpha=-4$ and
changing variables
$$
\{x={\frac {i}{t}},\ y={\frac {i\,t\,\left (t\,\u -1\right)}{{t}^{2}+1 }}\}.
$$

\section{Conclusions}
\label{conclusions}

In this paper, a multi-parameter non-constant-invariant Abel class has been
presented which generalizes in one the integrable cases shown in the works
by Abel, Liouville, Appell and others, including
all those shown in Kamke's book and other references. This new class
splits into various subclasses, many of which are fully integrable,
including some not previously shown elsewhere to the best of our
knowledge. 

The particular subclass represented by \eq{AIL_4} mapping 
non-constant-invariant Abel equations into linear first-order equations contains by itself most of
the exactly integrable cases we have seen. The subclass represented
by \eq{AIR_10} appears to us to be the most general class mapping Abel equations into 
Riccati ones; indeed it includes the parameterized mappings presented by
Liouville and Appell as particular members.

Finally, the mapping of Abel classes into other Abel classes
presents a useful way of finding new integrable classes from other classes
known to be solvable, as shown in \se{fitting}. We are presently working on
analyzing different connections between all these subclasses and expect to
find reportable results in the near future.

\ni {\bf Acknowledgments}

\noindent This work was supported by the Centre of Experimental and
Constructive Mathematics, Simon Fraser University, Canada, and by the
Symbolic Computation Group, Faculty of Mathematics, University of Waterloo,
Ontario, Canada. The authors would like to thank K. von
B\"ulow\footnote{Symbolic Computation Group of the Theoretical Physics
Department at UERJ - Brazil.} for a careful reading of this paper, and T.
Kolokolnikov\footnote{Department of Mathematics, University of British
Columbia - Canada.} and A. Wittkopf\footnote{Centre for Experimental and
Constructive Mathematics, Simon Fraser University - Canada.}
for useful related discussions.
We also acknowledge T. Kolokolnikov and
one anonymous referee for their suggestion of normalizing AIR by
normalizing the cubic polynomial on the numerator of \eq{AIL_8} and
\eq{AIR_10}.

\end{document}